\documentclass[preprint,oneside,12pt]{elsarticle}

\usepackage[utf8]{inputenc}
\usepackage{graphics}
\usepackage{graphicx}
\usepackage{epsfig}
\usepackage{latexsym}
\usepackage{amsthm}
\usepackage{amssymb}
\usepackage{amsmath}
\usepackage{amsfonts}
\usepackage{indentfirst}
\newcommand{\tmsamp}[1]{\textit{#1}}
\theoremstyle{definition}
 \newtheorem{Def}{Definition}
 \newtheorem{The}{Theorem}
 \newtheorem{Pro}{Proposition}
 \newtheorem{Lem}{Lemma}
 
 \newtheorem{Rem}{Remark}

 \newtheorem{Not}{Notation}

\DeclareMathOperator{\dd}{d}

\DeclareMathOperator{\L1}{L}
\DeclareMathOperator{\H1}{H}
\def\Rset{{\mathbb R}}

\journal{Ordinary Diffential Equation}

\begin{document}

\begin{frontmatter}

\title{Existence of solutions for a problem of resonance road space
with weight
\tnoteref{label1}}
\tnotetext[label1]{2000 Mathematics Subject Classification, 35J20}

\author{Antonio Ronaldo G.~Garcia, Moisés D.~dos Santos\fnref{label2}
and Adrião D.~D.~Neto\fnref{label3}}
\ead{ronaldogarcia@ufersa.edu.br, mdantas@ufersa.edu.br and adriao@dca.ufrn.br}
\ead[url]{www.ufersa.edu.br and www.dca.ufrn.br}

\address{Universidade Federal Rural do Semi-Árido -- UFERSA\fnref{label2}\\
Universidade Federal do Rio Grande do Norte -- UFRN\fnref{label3}}
\fntext[label2]{Departamento de Ciências Exatas e Naturais -- DCEN, 
CEP.: 59625-900, Mossoró-RN}
\fntext[label3]{Departamento de Engenharia de Computação e Automação
  -- DCA, CEP.: 59078-900, Natal-RN}

\begin{abstract}
In this work we use variational methods to show the existence of weak 
solutions for a nonlinear problem of the type elliptic. This problem
was initially study by the authors Ahmad, Lazer and Paul (see \cite{ALP})
considering the space $\Omega\subset\Rset^n$ a bounded domains. In this work
we extend your result now considering the domain $\Rset^n$. Indeed,
the main theorems in this paper constitute an extension to $\Rset^n$
of your previous results in bounded domains.                              
\end{abstract}


\end{frontmatter}


\section*{Introduction}\label{sec-1}



In this work we obtain a result of existence of weak solution
for the problem 
\begin{equation}\label{eq-1}
\left\{\begin{array}{ll}
-\Delta u+u=\lambda_k h(x)u+g(x,u), &\mbox{on}~\Rset^n,~n\ge 3\\
\qquad\quad ~u\in\H1^1(\Rset^n) &
\end{array}\right.
\end{equation}
where $h:\Rset^n\to\Rset_+,~g:\Rset^n\times\Rset\to\Rset$, are
continuous functions satisfying the following conditions:
\begin{enumerate}
\item[$(i)$] $h\in\L1^1(\Rset^n)\cap\L1^{\infty}(\Rset^n)$;
\item[$(ii)$] There is a function $Z\in
\L1^1(\Rset^n)\cap\L1^\infty(\Rset^n)$ such that $Z(x)>0,~\forall
~x\in\Rset^n$, and 
\begin{equation}\label{eq-2}
|g(x,t)|\le Z(x),~\forall~(x,t)\in\Rset^n\times\Rset.
\end{equation}
\end{enumerate}
Furthermore, $\lambda_k$ is the k-th eigenvalue associated with the
problem 
\begin{equation}\label{eq-3}
\left\{\begin{array}{ll}
-\Delta u+u=\lambda h(x)u, & \mbox{on}~\Rset^n, ~n\ge 3\\ 
\qquad\quad ~u\in\H1^1(\Rset^N)
\end{array}\right.
\end{equation}
characterizing the problem (\ref{eq-1}) as a problem of resonance (see
Definition 1.2 pg. 2 of \cite{GA}).

Furthermore, we assume that $g$ satisfies one of the conditions
$(g_2^+)$ or $(g_2^-)$, this is,
\begin{equation}\label{eq-4}
\int_{\Rset^N}G(x,v(x))\dd
x\underset{\|v\|\to\infty}{\longrightarrow}\pm\infty,
\end{equation}
where $v\in N_{\lambda_k}$ the eigenspace associated with the
eigenvalue $\lambda_k$ and $$G(x,\cdot)=\int_0^s g(x,\tau)\dd\tau$$ is
the primitive of the function $g(x,\cdot)$.

We find a weak solution for the problem (\ref{eq-1}) determining
critical points of the energy functional associated
$\Phi:\H1^1(\Rset^n)\to \Rset$ defined by
$$\Phi(u)=\frac{1}{2}\int_{\Rset^n}\left(|\Delta u|^2+|u|^2\right)\dd
x-\frac{\lambda_k}{2}\int_{\Rset^n} u^2h\dd x-\int_{\Rset^n}G(x,u)\dd
x.$$

In this work we generalized the result from Lazer-Ahmad-Paul in 
\cite{ALP}. For this was necessary some theorics results once again that
the domain considered here is no bounded. In a space no bounded no
exists immersion compact of the Soboleve space $\H1^1(\Rset^n)$ in
$\L1^p(\Rset^n)$. We contorted this working in spaces with weight. 

Problems at resonance have been of interest to researchers ever since
the pioneering work of Landesman and Lazer \cite{LL} in 1970 for
second order elliptic operators in bounded domains. The literature on
resonance problems in bounded domains is quite vast; of particular
interest to this paper are the works of Ahmad, Lazer and Paul
\cite{ALP} in 1976 and of Rabinowitz in 1978, in which critical point
methods are applied. Recently, using other tecnic, Garza and 
Rumbos \cite{GA} do a result that is a extension to $\Rset^n$ 
of the Ahmad, Lazer and Paul result. In this paper, furthermore of
resonance they obtain result on strong resonance. 

Resonance problems on unbounded domains, and in particular in
$\Rset^n$, have been studied recently by Costa and Tehrani \cite{TC}
and by Jeanjean \cite{JJ}, and by Stuart and Zhou \cite{SZ} for
radially symmetric solutions for asymptotically linear problems in
$\Rset^n$. In all these references variational methods were
used. Hetzer and Landesman \cite{HL} for resonant problems for a class
of operators which includes the Schrödinger operator. The problem of
resonance that we work can be see in \cite{WIL} for bounded domain.      

We write in general the direction of our proposal, In the
Section \label{sec-2} we define spaces with weight and show that this
space is Banach. In the Section \ref{sec-3} we show two important
results, a of immersions continuous and the other of immersion
compact space of the space $\H1^1(\Rset^n)$ in $\L1^p(\Rset^n,h\dd
x)$. Finally, in the Section \ref{sec-4} we show the main result of 
this work that guarantees the existence of weak solution to the 
problem (\ref{eq-1}).

\section{Spaces with weight}\label{sec-2}
In this section we define spaces with weight and show some properties
involving these spaces. For this, start with the following:
\begin{Def}\label{prima-1} \tmsamp{Let $h:\Rset^N\to (0,+\infty)$ be 
a measurable function and $1<p<\infty$. We define the space 
$\L1^p(\Rset^N,hdx)$ the space of all
the measurable functions $f:\Rset^n\to\Rset$ such that
$\int_{\Rset^n}|f(x)|^ph(x)\dd x<\infty$, i.e, 
$$\L1^p\left(\Rset^N,hdx\right):=\left\{f:\Rset^N\to\Rset ~\mbox{measurable}| ~
\int_{\Rset^N}|f(x)|^ph(x)\dd x<\infty\right\}.$$}
\end{Def}

In the that follow we denoted by 
$$\left\|u\right\|=\left(\int_{\Rset^N}\left(\left|\nabla{u}\right|^2+
\left|u\right|^2\right)\dd x\right)^{{1}/{2}}~
\mbox{and}~\left\|f\right\|_{p,h}=\left(\int_{\Rset^N}\left|f\right|^ph
  \dd x
\right)^{1/p}$$
the norms in $\H1^1(\Rset^N)$ and $\L1^p(\Rset^N,h\dd x),$ respectively.
\begin{The}\label{prima-2} \tmsamp{The space
    $\left(\L1^p(\Rset^N,h\dd x);\|\cdot\|_{p,h}\right)$ with
    $1\leq {p}<\infty$ is a Banach's space.}
\end{The}
\begin{proof}
Let $\{u_{n}\}\subset{\L1^p(\Rset^N,hdx)}$ be a Cauchy's sequence.
Then, give $\epsilon>0$ there exists $n_0 \in{\mathbb{N}}$ such that
\begin{equation}\label{cauchy-1}
\|u_n-u_m\|_{p,h}<\epsilon, ~\forall~n,m\geq n_0.
\end{equation}
Define $v_{n}=h^{\frac{1}{p}}u_{n},$ it follow from (\ref{cauchy-1}) that
$$\|v_{n}-v_{m}\|_{\L1^p(\Rset^N)}<\epsilon, ~\forall ~ n,m\geq n_0.$$
Since $\L1^p(\Rset^N)$ is a Banach space, there exists
$w \in{\L1^p(\Rset^N)}$ with
\begin{equation}\label{cauchy-2}
v_{n}\underset{n\longrightarrow\infty}{\longrightarrow} w 
\end{equation}
Defining, $u(x)=\frac{w(x)}{h(x)^{\frac{1}{p}}},$
note that
\begin{eqnarray*}
\|u_{n}-u\|_{p,h}&=&\left(\int_{\Rset^N}\left|u_{n}-u\right|^{p}hdx\right)^{\frac{1}{p}}\\
&=&\left(\int_{\Rset^N}\left|h^{\frac{1}{p}}u_{n}-h^{\frac{1}{p}}u\right|^{p}dx\right)^{\frac{1}{p}}
\end{eqnarray*}
which implies,
\begin{eqnarray*}
\|u_{n}-u\|_{p,h}&=&\left(\int_{\Rset^N}\left|v_{n}-w\right|^{p}dx\right)^{\frac{1}{p}}\\
&=&\|v_{n}-w\|_{\L1^p(\Rset^N)}.
\end{eqnarray*}
From (\ref{cauchy-2}), we have
$\|u_{n}-u\|_{p,h}\underset{n\longrightarrow\infty}{\longrightarrow} 0$,
that is, $u_{n}\underset{n\longrightarrow\infty}{\longrightarrow} u$ in 
$\L1^p(\Rset^N,hdx)$. Therefore, $\L1^p(\Rset^N,hdx)$ is a
Banach's space.
\end{proof}

From Theorem \ref{prima-2} we can conclude that 
$\left(L^2(\Rset^N,hdx),\|\cdot\|_{2,h}\right)$
is a Hilbert's space with the inner product
$\left(f,g\right)_{2,h}=\int_{\Rset^N}hfg\dd x$.
\section{Auxiliary results}\label{sec-3}
In this section we show two auxiliary results important for we
show that the problem (\ref{eq-1}) have weak solution. 
The first is a result of continuous immersion and the second a result 
of compact immersion,  
both of the space $\H1^1(\Rset^n)$ in $\L1^p(\Rset^n,h\dd x)$.
\subsection{A result of continuous immersion} 
\begin{The}\label{prima-3} 
\tmsamp{If $h \in L^{\infty}\left(\Rset^N\right)$,
    then applies continuous immersion 
$$\H1^1(\Rset^N)\hookrightarrow \L1^p(\Rset^N,hdx)$$
for all $p\in \left[1,2^*\right]~(2^*=2N/(N-2))$ if $N \geq 3$.}
\end{The}

\begin{proof} For $u \in \H1^1(\Rset^N)$, we have
\begin{eqnarray*}
\left(\int_{\Rset^N}\left|u\right|^phdx\right)^{\frac{1}{p}}&\leq&\left\|h\right\|_{\infty}\left(\int_{\Rset^N}\left|u\right|^p\dd
  x\right)^{\frac{1}{p}}\\
&=&\left\|h\right\|_{\infty}\left\|u\right\|_{\L1^p(\Rset^N)}.
\end{eqnarray*}
Since $h \in L^{\infty}(\Rset^N)$. From Sobolev's continuous
immersion, we have $$\H1^1(\Rset^N)\hookrightarrow
\L1^p(\Rset^N), ~\forall~p\in \left[2,2^*\right],$$
if $N \geq 3$. Thus, there exists $C>0$ such that
$$\left\|u\right\|_{\L1^p(\Rset^N)}\leq C\left\|u\right\|,~\forall
~u\in \H1^1(\Rset^N).$$
Therefore,
$$\left(\int_{\Rset^N}\left|u\right|^phdx\right)^{1/p}\leq
C\left\|h\right\|_{\infty}\left\|u\right\|, ~\forall~ u\in \H1^1(\Rset^N).$$
Considering $C_{1}=C\left\|h\right\|_{\infty}$, we have
$$\left(\int_{\Rset^N}\left|u\right|^phdx\right)^{1/p}\leq C_{1}
\left\|u\right\|,$$
that is,
$$\left\|u\right\|_{p,h}\leq C_{1}\left\|u\right\|, ~\forall ~u\in
\H1^1(\Rset^N)$$
showing the continuous immersion. 
\end{proof}
\subsection{A result of Compact immersion}
\begin{Lem}\label{prima_4.1} \tmsamp{Give $\epsilon>0$, there exists
    $R>0$ such that 
$$\|u_{n_{j}}-u\|_{\L1^p(B^{c}_{R}(0),hdx)}<\frac{\epsilon}{2}, ~\forall~ n_{j}.$$}
\end{Lem}
\begin{proof} Indeed, to $R>0 $ and using the Hölder inequality with
exponents
$\frac{1}{\gamma}+ \frac{1}{\frac{2^*}{p}}=1$, we get
$$
\int_{B^{c}_{R}(0)}\left|u_{n_{j}}-u\right|^{p}h\dd x\leq \left\|h\right\|_{L^{\gamma}(B^{c}_{R}(0))}.\left\|\left|u_{n_{j}}-u\right|^{p}\right\|_{L^{\frac{2*}{p}(B^{c}_{R}(0))}}
$$
which implies,
$$
\int_{B^{c}_{R}(0)}\left|u_{n_{j}}-u\right|^{p}h\dd x\leq \left\|h\right\|_{L^{\gamma}(B^{c}_{R}(0))
}.\left\|u_{n_{j}}-u\right\|^{p}_{L^{{2}^{*}}(\Rset^N)}
$$
and from the Sobolev's continuous immersions, we have
$$
\int_{B^{c}_{R}(0)}\left|u_{n_{j}}-u\right|^{p}h\dd x\leq C\left\|h\right\|_{L^{\gamma}(B^{c}_{R}(0))}.\left\|u_{n_{j}}-u\right\|^{p}.
$$
Consequently,
\begin{equation}\label{compac-1}
\int_{B^{c}_{R}(0)}h\left|u_{n_{j}}-u\right|^{p}\dd x\leq C_{1}\left\|h\right\|_{L^{\gamma}(B^{c}_{R}(0))}.
\end{equation}
Follows the theory of measure that if $h\in
L^{\gamma}(\Rset^{N}), ~\forall~ \gamma \in [1,\infty)$, give
$\epsilon>0$ there exists $R>0$ such that 
$$
\left\|h\right\|_{L^{\gamma}(B^{c}_{R}(0))}<\left(\frac{\epsilon}{2}\right)^{p}
\frac{1}{C_{1}}.
$$
Now, from result above and from (\ref{compac-1}), we obtain
$$
\int_{B^{c}_{R}(0)}h\left|u_{n_{j}}-u\right|^{p}dx\leq \left(\frac{\epsilon}{2}\right)^{p}\;\; \forall\;\; n_{j}.
$$
Therefore,
\begin{equation}\label{compac-2}
\left\|u_{n_{j}}-u\right\|_{\L1^p(B^{c}_{R}(0),h\dd
  x)}<\frac{\epsilon}{2}~\forall ~n_{j}
\end{equation}
which show the lemma.
\end{proof}

\begin{The}\label{prima-4} \tmsamp{If $h \in L^{1}({\Rset^N})\cap
    L^{\infty}({\Rset^N}),$ has the compact immersion
$$\H1^{1}({\Rset^N})\hookrightarrow  \L1^p(\Rset^N,h\dd x),
~\forall~p \in [1,2^{*})$$ if $N\geq 3$.}
\end{The}
\begin{proof}
Let $(u_{n})\subset \H1^{1}({\Rset^N})$ be a bounded sequence,
using the fact that $\H1^{1}({\Rset^N})$ is a reflexive space
there exists $(u_{n_{j}})\subset (u_{n})$ such that
$$u_{n_{j}}\rightharpoonup u ~\mbox{in}~ \H1^{1}({\Rset^N}).$$
Furthermore,
\begin{eqnarray*}
\|u_{n_{j}}-u\|_{\L1^p(B_{R}(0),h\dd
  x)}&=&\left(\int_{B_{R}(0)}h\left|u_{n_{j}}-u\right|^{p}\dd x\right)^
{\frac{1}{p}}\\
&\leq& \|h\|_{\infty}.\|u_{n_{j}}-u\|_{\L1^p(B_{R}(0))}.
\end{eqnarray*}
From Sobolev's compact immersions, we have
$\H1^1(\Rset^N)\hookrightarrow  \L1^p(B_{R}(0))$, of where follows
$$u_{n_{j}}\underset{n_j\longrightarrow\infty}{\longrightarrow} u~\mbox{in}~ \L1^p(B_{R}(0)).$$
Thus, give $\epsilon >0 $ there exists $n_{j_{0}}$ such that
$$\|u_{n_{j}}-u\|_{\L1^p(B_{R}(0))}<\frac{\epsilon}{2\|h\|_{\infty}}, ~\forall~ n_{j}\geq n_{j_{0}},$$
thus, we have
\begin{equation}\label{compac-3}
\|u_{n_{j}}-u\|_{\L1^p(B_{R}(0),h\dd x)}<\frac{\epsilon}{2}, ~\forall~
 n_{j}\geq n_{j_{0}}.
\end{equation}
From Lemma \ref{prima_4.1} and (\ref{compac-3}), we have
$$
\left\|u_{n_{j}}-u\right\|_{p,h}<\epsilon, ~\forall~ n_{j}\geq n_{j_{0}}
$$
which implies,
$$
u_{n_{j}}\underset{n_j\longrightarrow\infty}{\longrightarrow} u~
\mbox{in}~ \L1^p(\Rset^N,h\dd x)$$
showing the compactness.
\end{proof}
\section{Main result}\label{sec-4}

\subsection{Preliminaries} 
We consider the functional 
$\Phi :\H1^1(\Rset^N)\to \Rset$ such that 
\begin{equation}\label{energia-1}
\Phi(u)=\frac{1}{2}\int_{\Rset^N}(|\nabla 
u|^{2}+|u|^{2})\dd x-\frac{\lambda_{k}}{2}
\int_{\Rset^N}u^{2}h\dd x-\int_{\Rset^N}G(x,u)\dd x.
\end{equation}
This functional is well defined in $\H1^1(\Rset^N)$, is of class   
$C^{1}(\H1^1(\Rset^N),\Rset)$ (see \cite{MOI}) 
and its critical points are weak solutions of (\ref{eq-1}). 

We set the linear operator $Lu=u-\lambda_{k}S(u)$, we have
$$(Lu,u)_{\H1^1(\Rset^N)}=\int_{\Rset^N}(\nabla(Lu)\nabla u
+Lu.u)\dd x,$$
thus,
$$(Lu,u)_{\H1^1(\Rset^N)}=\int_{\Rset^N}(\nabla(u-\lambda_{k}S(u))\nabla
u\dd x +\int_{\Rset^N}(u-\lambda_{k}S(u))u\dd x,$$
this is,
\begin{equation}\label{energia-2}
(Lu,u)_{\H1^1(\Rset^N)}=\int_{\Rset^N}(|\nabla
u|^{2}+|u|^{2})\dd x - \lambda_{k}\int_{\Rset^N}(\nabla
S(u)\nabla u+S(u)u)\dd x.
\end{equation}
Moreover, considering $S(u)=w$, we have 
$$\left\{\begin{array}{ll}
-\Delta{w}+w=hu,& \Rset^N\\
\qquad\quad ~~w \in \H1^1(\Rset^N) &
\end{array}\right.
$$
which implies,
$$\int_{\Rset^N}(\nabla\phi\nabla w+\phi w)\dd
x=\int_{\Rset^N}uh\phi\dd x, ~\forall~\phi \in \H1^1(\Rset^N).$$
Fixing $u=\phi$,
\begin{equation}\label{energia-3}
\int_{\Rset^N}(\nabla S(u) \nabla u +S(u)u)\dd
x=\int_{\Rset^N}u^{2}h\dd x.
\end{equation}
Therefore, from (\ref{energia-2}) and (\ref{energia-3}), we have
$$
\frac{1}{2}(Lu,u)_{\H1^1(\Rset^N)}=\frac{1}{2}\int_{\Rset^N}(|\nabla
u|^{2} + |u|^{2})\dd x
-\frac{\lambda_{k}}{2}\int_{\Rset^N}u^{2}h\dd x
$$
from which follows
$$\Phi(u)=\frac{1}{2}(Lu,u)_{\H1^1(\Rset^N)}-\int_{\Rset^N}G(x,u)\dd
x.$$

Now, we fix the following orthogonal decomposition of $X=\H1^1(\Rset^N)$,
$$X=X_{-}\oplus X_{0}\oplus X_{+},$$ where $$X_{0}=N_{\lambda_{k}},~
X_{-}=N_{\lambda_{1}}\oplus N_{\lambda_{2}}\oplus...\oplus N_{\lambda_{k-1}}
~\mbox{and} ~X_{+}=N_{\lambda_{k+1}}\oplus N_{\lambda_{k+2}}\oplus\dots$$

\begin{Pro}\label{main-1} \tmsamp{If $u \in X_{0}$, then 
$(Lu,u)_{\H1^1\left(\Rset^N\right)}=0$.}
\end{Pro}
\begin{proof} As we have seen that
$$
(Lu,u)_{\H1^1(\Rset^N)}=\int_{\Rset^N}(|\nabla
u|^{2}+|u|^{2})\dd x-\lambda_{k}\int_{\Rset^N}u^{2}h\dd x,$$
since $u \in X_{0},$ it is solution of the problem
$$
 \left\{\begin{array}{rcl}
-\Delta{u}+u &\mbox{=}&\lambda_{k}hu,~\Rset^N \\
u            &\in & \H1^1(\Rset^N)
\end{array}\right.
$$
of where, we have
$$\int_{\Rset^N}(\nabla\phi\nabla u +\phi
u)=\lambda_{k}\int_{\Rset^N}u\phi h\dd x, ~\forall~\phi\in \H1^1\left(\Rset^N\right).$$
Fixing $\phi=u,$
$$\int_{\Rset^N}(|\nabla u|^{2}+|u|^{2})\dd x=
\lambda_{k}\int_{\Rset^N}u^{2}h\dd x,$$
showing that $$(Lu,u)_{\H1^1(\Rset^N)}=0,$$ 
if $u\in X_{0}$.
\end{proof}

\begin{Pro}\label{main-2}\tmsamp{If $u\in X_{-}$, then there exists
$\alpha>0$ such that $$(Lu,u)_{\H1^1(\Rset^N)}\leq-\alpha\|u\|^{2},$$ 
this is,  $L$ is negative defined in $X_{-}$.}
\end{Pro}
\begin{proof}
Let $u\in X_{-}=N_{\lambda_{1}}\oplus N_{\lambda_{2}}\oplus...\oplus
N_{\lambda_{k-1}},$ be, then 
$$
u=\phi_{1}+\phi_{2}+...+\phi_{k-1}~\mbox{and}~\nabla u=\nabla \phi_{1}+ \nabla \phi_{2}+...+ \nabla \phi_{k-1}.
$$
Note that from (\ref{energia-2})
$$(Lu,u)_{\H1^1(\Rset^N)}=\int_{\Rset^N}(|\nabla
u|^{2}+|u|^{2})\dd x-\lambda_{k}\int_{\Rset^N}(\nabla S(u)\nabla
u+uS(u))\dd x$$
thus,
\begin{equation}\label{energia-vio}
(Lu,u)_{\H1^1(\Rset^N)}=(u,u)_{\H1^1(\Rset^N)}-\lambda_{k}(S(u),u)_{\H1^1(\Rset^N)}.
\end{equation}
Provided that $\phi_{j}$ satisfies
$$
 \left\{
\begin{array}{rcl}
-\Delta\phi_{j}+\phi_{j} &\mbox{=}&\lambda_{j}h\phi_{j} ,\; \Rset^N \\
\phi_{j} & \in & \H1^1(\Rset^N)
\end{array}
\right.
$$
From definition of the solution operator
$$
S(\lambda_{j}\phi_{j})=\phi_{j}
$$
therefore,
$$S(\phi_{j})=\frac{1}{\lambda_{j}}\phi_{j}~\mbox{and}~\nabla 
S(\phi_{j})=\frac{1}{\lambda_{j}}\nabla\phi_{j}.$$
From linearity of $S$, we have
$$S(u)=S(\phi_{1})+S(\phi_{2})+\dots+S(\phi_{k-1})$$
which implies,
$$S(u)=\frac{1}{\lambda_{1}}\phi_{1}+\frac{1}{\lambda_{2}}\phi_{2}+\dots+\frac{1}{\lambda_{k-1}}\phi_{k-1}.$$
Thus,
$$(S(u),u)_{\H1^1(\Rset^N)}=(S(\phi_{1})+S(\phi_{2})+\dots+S(\phi_{k-1}),\phi_{1}+\phi_{2}+\dots+\phi_{k-1})_{\H1^1(\Rset^N)}$$
and provided that $(\phi_{j},\phi_{k})_{\H1^1(\Rset^N)}=0$, if
$j\neq k$, we obtain
$$(S(u),u)_{\H1^1(\Rset^N)}=\frac{1}{\lambda_{1}}(\phi_{1},\phi_{1})_{\H1^1(\Rset^N)}+\dots+\frac{1}{\lambda_{k-1}}(\phi_{k-1},\phi_{k-1})_{\H1^1(\Rset^N)}.$$
It follows from definition of inner product in $\H1^1(\Rset^N)$ that
$$
(S(u),u)_{\H1^1(\Rset^N)}=\frac{1}{\lambda_{1}}\int_{\Rset^N}(|\nabla\phi_{1}
|^{2}+|\phi_{1}|^{2} )\dd
x+\dots+\frac{1}{\lambda_{k-1}}\int_{\Rset^N}(|\nabla\phi_{k-1}|^{2}+|\phi_{k-1}|^{2})\dd
x$$
which implies
$$\lambda_{k}(S(u),u)_{\H1^1(\Rset^N)}=\sum^{k-1}_{j=1}\int_{\Rset^N}\frac{\lambda_{k}}{\lambda_{j}}(|\nabla\phi_{j}|^{2}+|\phi_{j}|^{2})\dd
x,$$
this is,
$$
\lambda_{k}\left(S(u),u\right)_{\H1^1(\Rset^N)}=\sum^{k-1}_{j=1}\frac{\lambda_{k}}{\lambda_{j}}\left\|\phi_{j}\right\|^{2}.
$$
Furthermore,
$$(u,u)_{\H1^1(\Rset^N)}=(\phi_{1}+\phi_{2}+\dots+\phi_{k-1},\phi_{1}+\phi_{2}+\dots+\phi_{k-1})_{\H1^1(\Rset^N)}$$
hence,
$$(u,u)_{\H1^1(\Rset^N)}=\sum^{k-1}_{j=1}\int_{\Rset^N}(|\nabla\phi_{j}|^{2}
+|\phi_{j}|^{2})\dd x=\sum^{k-1}_{j=1}\|\phi_{j}\|^{2}.$$
Now, we can to writer (\ref{energia-vio}) of the follows way
$$(Lu,u)_{\H1^1(\Rset^N)}=\sum^{k-1}_{j=1}\left(1-\frac{\lambda_{k}}{\lambda_{j}}\right)\|\phi_{j}\|^{2}$$
and using the fact that
$$
\lambda_{j}\leq\lambda_{k-1}, ~\forall~ j \in \{1,2,...,k-1\}
$$
which implies,
$$1-\frac{\lambda_{k}}{\lambda_{k-1}}\geq 1-\frac{\lambda_{k}}{\lambda_{j}}$$
thus,
\begin{eqnarray*}
(Lu,u)_{\H1^1(\Rset^N)}&\leq&\sum^{k-1}_{j=1}\left(1-\frac{\lambda_{k}}{\lambda_{k-1}}\right)\|\phi_{j}\|^{2}\\
&=&\left(1-\frac{\lambda_{k}}{\lambda_{k-1}}\right)\sum^{k-1}_{j=1}\|\phi_{j}\|^{2}.
\end{eqnarray*}
Therefore,
$$(Lu,u)_{\H1^1(\Rset^N)}\leq\left(1-\frac{\lambda_{k}}{\lambda_{k-1}}\right)
\|u\|^{2}$$
from where we have
$$(Lu,u)_{\H1^1(\Rset^N)}\leq-\left(\frac{\lambda_{k}}{\lambda_{k-1}}-1\right)
\|u\|^{2},$$
this is,
$$(Lu,u)_{\H1^1(\Rset^N)}\leq-\alpha\|u\|^{2}<0, ~\forall~ u \in
X_{-}\setminus \{0\}$$ and $\alpha>0$.
\end{proof}

\begin{Pro}\label{main-3} \tmsamp{If $u\in X_{+}$, then there exists
    $\alpha>0$ such that
    $(Lu,u)_{\H1^1(\Rset^N)}\geq\alpha\|u\|^{2}$, this is, $L$ is
    positive defined in $X_{+}$.}
\end{Pro}

\begin{proof}
Suppose that $u\in X_{+}=N_{\lambda_{k+1}}\oplus
N_{\lambda_{k+2}}\oplus\dots$. Then $u=\phi_{k+1}+\phi_{k+2}+\dots$
We consider, $u=\sum^{\infty}_{j=k+1}\phi_{j}=\lim_{N \to
  \infty}w_{N}$, where $w_{N}=\sum^{N}_{j=k+1}\phi_{j}, ~w_{N}\in X_{+}$.
We show that there exists $\alpha>0$ independent of $N$ such that
$$(Lw_{N}, w_{N})_{\H1^1(\Rset^N)}\geq\alpha\|w_{N}\|^{2},
~\forall~ N\in \{k+1,k+2,\dots\}.$$
Indeed, note that
$$
(Lw_{N},w_{N})_{\H1^1(\Rset^N)}=\left(\sum^{N}_{j=k+1}L\phi_{j},\sum^{N}_{j=k+1}\phi_{j}\right)_{\H1^1(\Rset^N)}
$$
which implies,
$$(Lw_{N},w_{N})_{\H1^1(\Rset^N)}=\left(\sum^{N}_{j=k+1}(\phi_{j}-\lambda_{k}S(\phi_{j})),\sum^{N}_{j=k+1}\phi_{j}\right)_{\H1^1(\Rset^N)}$$
consequently,
$$(Lw_{N}, w_{N})_{\H1^1(\Rset^N)}=\left(\sum^{N}_{j=k+1}\phi_{j},\sum^{N}_{j=k+1}\phi_{j}\right)_{\H1^1(\Rset^N)}-\left(\sum^{N}_{j=k+1}\lambda_{k}S(\phi_{j}),\sum^{N}_{j=k+1}\phi_{j}\right)_{\H1^1(\Rset^N)}.$$
Now, using the same reasoning from Proposition \ref{main-2}, we have
$$(Lw_{N},w_{N})_{\H1^1(\Rset^N)}=\sum^{N}_{j=k+1}\|\phi_{j}\|^{2}-\sum^{N}_{j=k+1}\frac{\lambda_{k}}{\lambda_{j}}\|\phi_{j}\|^{2}$$
and noting that $\lambda_{k}\geq\lambda_{k+1}$ we find
$$1-\frac{\lambda_{k}}{\lambda_{j}}\geq 1-\frac{\lambda_{k}}{\lambda_{k+1}}.$$
Therefore,
$$(Lw_{N},w_{N})_{\H1^1(\Rset^N)}\geq\left(1-\frac{\lambda_{k}}{\lambda_{k+1}}\right)\sum^{N}_{j=k+1}\|\phi_{j}\|^{2}.$$
Thus,
$$(Lw_{N},w_{N})_{\H1^1(\Rset^N)}\geq\left(1-\frac{\lambda_{k}}{\lambda_{k+1}}\right)\|w_{N}\|^{2}$$
showing that
\begin{equation}\label{energia-4}
(Lw_{N},w_{N})_{\H1^1(\Rset^N)}\geq\alpha\|w_{N}\|^{2}, ~\forall ~N
\in \{k+1,k+2,\dots\},
\end{equation}
where $\alpha=\left(1-\frac{\lambda_{k}}{\lambda_{k+1}}\right)$.
Using the fact that
$w_{N}\underset{N\rightarrow\infty}{\longrightarrow} u$, we obtain
$$\|w_{N}\|^{2}\underset{N\rightarrow\infty}{\longrightarrow} \|u\|^{2}$$
and from continuity of $L$, we obtain
$L(w_{N})\underset{N\rightarrow\infty}{\longrightarrow} L(u)$.
Consequently,
$$(Lw_{N},w_{N})_{\H1^1(\Rset^N)}\underset{N\rightarrow\infty}{\longrightarrow}(Lu,u)_{\H1^1(\Rset^N)}$$
thus, passing to the limit at (\ref{energia-4}) follows that
$(Lu,u)_{\H1^1(\Rset^N)}\geq\alpha\|u\|^{2}$
\end{proof}

\begin{Not}\label{main-5} \tmsamp{In the that follows, we denoted by
$P_{-},P_{0}$ and $P_{+}$ the orthogonais projections on $X_{-},X_{0}$
and $X_{+}$, respectively.}
\end{Not}

\begin{Lem}\label{main-6} \tmsamp{If $u=P_{0}u+P_{+}u \in X_{0}\oplus
    X_{+}$, then 
$(Lu,u)_{\H1^1(\Rset^N)}=(LP_{+}u,P_{+}u)_{\H1^1(\Rset^N)}$.}
\end{Lem}
\begin{proof} Suppose that
$(Lu,u)_{\H1^1(\Rset^N)}=(L(P_{0}u+P_{+}u),P_{0}u+P_{+}u)_{\H1^1(\Rset^N)}$. 
Then  
$$(Lu,u)_{\H1^1(\Rset^N)}=(LP_{0}u+LP_{+}u,P_{0}u+P_{+}u)_{\H1^1(\Rset^N)}.$$
Using the Proposition \ref{main-1} and the fact that $L$ is symmetric
operator, we obtain
$$(Lu,u)_{\H1^1(\Rset^N)}=2(LP_{0}u,P_{+}u)_{\H1^1(\Rset^N)}+(LP_{+}u,P_{+}u)_{\H1^1(\Rset^N)}.$$
Note that
$(LP_{0}u,P_{+}u)_{\H1^1(\Rset^N)}=(P_{0}u-\lambda_{k}S(P_{0}u),P_{+}u)_{\H1^1(\Rset^N)}$
hence,
$$(LP_{0}u,P_{+}u)_{\H1^1(\Rset^N)}=\lambda_{k}(S(P_{0}u),P_{+}u)_{\H1^1(\Rset^N)}.$$
Therefore, by the definition of operator solution, we obtain
$$(LP_{0}u,P_{+}u)_{\H1^1(\Rset^N)}=(P_{0}u,P_{+}u)_{L^2(\Rset^N,hdx)}=0,$$
where the last equality we use the orthogonality of the projections. Thus,
$(Lu,u)_{\H1^1(\Rset^N)}=(LP_{+}u,P_{+}u)_{\H1^1(\Rset^N)}.$
\end{proof}

\begin{Pro}\label{main-7} \tmsamp{Suppose valid the conditions
    (\ref{eq-2}) and $(g^{-}_{2})$. Then
\begin{enumerate}
\item[$(a)$]
$\Phi(u)\underset{\|u\|\rightarrow +\infty}{\longrightarrow} -\infty, ~u \in X_{-}$.
\item[$(b)$] $\Phi(u)\underset{\|u\|\rightarrow +\infty}{\longrightarrow} +\infty,~u \in X_{0}\oplus X_{+}$.
\end{enumerate}}
\end{Pro}
\begin{proof}
$(a)$ Suppose that $u\in X_{-}$. Using the fact that $L$ is defined
negative in $X_{-}$, we obtain
$$\Phi(u)\leq
-\frac{1}{2}\alpha\|u\|^{2}-\int_{\Rset^{N}}G(x,u)\dd x.$$
Note that for the mean value Theorem, we have
$|G(x,u)-G(x,0)|=|G^{'}(x,s)||u|$, $s(x)\in[0,u(x)]$ or $s(x)\in[u(x),0]$,
thus from (\ref{eq-2})
$$|G(x,u)-G(x,0)|\leq Z(x)|u|$$ therefore,
$-\int_{\Rset^{N}}G(x,u)\dd x\leq \|u\|_{1,Z}$,
of where we have that
$$\Phi(u)\leq -\frac{1}{2}\alpha\|u\|^{2}+\|u\|_{1,Z}.$$
From Theorem \ref{prima-4} follows that
$\|u\|_{2,Z}\leq C_{1}\|u\|$, then
$$\Phi(u)\leq -\frac{1}{2}\alpha\|u\|^{2}+C_{2}\|u\|$$
consequently,
$\Phi(u)\underset{\|u\|\rightarrow +\infty}{\longrightarrow} -\infty,
~u 
\in X_{-}$.

$(b)$ Suppose that $u=P_{0}u+P_{+}u \in X_{0}\oplus X_{+}$. Using the
Lemma \ref{main-6} has been
$\Phi(u)=\frac{1}{2}(LP_{+}u,P_{+}u)_{\H1^1(\Rset^N)}-\int_{\Rset^{N}}
G(x,u)\dd x$.
Furthermore, since $L$ defined positive in $X_{+}$ follows that
$$\Phi(u)\geq\frac{1}{2}\alpha\|P_{+}u\|^{2}-\int_{\Rset^{N}}
[G(x,u)-G(x,P_{0}u)]\dd x-\int_{\Rset^{N}}G(x,P_{0}u)\dd x.$$
Using again the mean value Theorem and (\ref{eq-2}), we have
\begin{eqnarray*}
-\int_{\Rset^{N}}[G(x,u)-G(x,P_{0}u)]\dd x&\leq& 
\int_{\Rset^{N}}|G(x,u)-G(x,P_{0}u)|\dd x\\
&\leq&\int_{\Rset^{N}}|u-P_{0}u|Z(x)\dd x
\end{eqnarray*}
which implies,
$$-\int_{\Rset^{N}}[G(x,u)-G(x,P_{0}u)\dd x\leq\|P_{+}u\|_{1,Z}.$$
Therefore, follows from the continuous immersion (see Theorem
\ref{prima-3}) that
\begin{equation}\label{energia-5}
\Phi(u)\geq\frac{1}{2}\alpha\|P_{+}u\|^{2}-C\|P_{+}u\|-\int_{\Rset^{N}}
G(x,P_{0}u)\dd x.
\end{equation}
Now, using the condition $(g^{-}_{2})$, the fact that
$\|u\|^{2}=\|P_{+}u\|^{2}+\|P_{0}u\|^{2}$
and analyzing the cases:
\begin{enumerate}
\item[$(i)$] $\|P_{0}u\|\rightarrow +\infty$ and $\|P_{+}u\|\leq M$. 
Using $(g^{-}_{2})$ and doing the analysis in (\ref{energia-5}), we have
$\Phi(u)\underset{\|u\|\rightarrow +\infty}{\longrightarrow} +\infty$.
\item[$(ii)$] $\|P_{+}u\|\rightarrow +\infty$ e $\|P_{0}u\|\leq K$. Note that
\begin{eqnarray*}
\left|-\int_{\Rset^{N}}G(x,P_{0}u)\dd x\right|&\leq&\int_{\Rset^{N}}|G(x,P_{0}u)|\dd x\\
&\leq&\int_{\Rset^{N}}|P_{0}u|Z(x)\dd x,
\end{eqnarray*}
that is,
$$\left|-\int_{\Rset^{N}}G(x,P_{0}u)dx\right|\leq K_{1}, ~K_1\in\Rset_+.$$
Thus, doing the analysis in (\ref{energia-5}), we have
$\Phi(u)\underset{\|u\|\rightarrow +\infty}{\longrightarrow} +\infty$.
\item[$(iii)$] If $\|P_{+}u\|\rightarrow +\infty$ and $\|P_{0}u\|\rightarrow +\infty$. Again the condition $(g^{-}_{2})$ and doing the analysis in (\ref{energia-5}, follows that
$\Phi(u)\underset{\|u\|\rightarrow +\infty}{\longrightarrow} +\infty$,
\end{enumerate}
which shows $(b)$. 
\end{proof}

\begin{Rem}\label{main-8} \tmsamp{If we had assumed the condition 
$(g^{+}_{2})$ to $(g^{-}_{2})$ instead of the reasoning used was the 
same and the conclusions were the following:
\begin{enumerate}
\item[$(a)$] $\Phi(u)\underset{\|u\|\rightarrow +\infty}{\longrightarrow} 
-\infty, ~u \in X_{0}\oplus X_{-}$.
\item[$(b)$] $\Phi(u)\underset{\|u\|\rightarrow +\infty}\longrightarrow 
+\infty, ~u \in X_{+}$.
\end{enumerate}}
\end{Rem}

\subsection{Weak solution to the problem (\ref{eq-1})}
We are ready to state and prove the main result of our work. This
result guarantees the existence of a critical point to the functional
$\Phi$, and therefore the existence of a weak solution to the problem
(\ref{eq-1}). 


For we show the existence of weak solution for the problem
(\ref{eq-1}) we use the theorem from the saddle of Rabionowitz.
\begin{The}[Saddle Point Theorem \cite{PHR}]
\tmsamp{Let $X=V\oplus W$ be a Banach' space, of way that 
$dim\ V<\infty,$ and let $\phi \in C^{1}(X,\Rset)$ be a  
maps satisfying the condition of Palais-Smale. 
If $D$ is a bounded neighborhood of the $0$ in $V$ such that
\begin{equation}\label{rabionowwitz}
a=\max_{\partial D}\phi<\inf_{W} \phi\equiv b,
\end{equation}
then
$$c=\inf_{h\in \Gamma }\max_{u\in \overline{D}}\phi(h(u))$$
is a critical value of $\phi$ with $c\geq b$, where
$$\Gamma=\left\{h\in C\left(\overline{D},X\right)|h(u)=u, 
~\forall~u\in\partial D\right\}.$$}
\end{The}

\begin{The}\label{main} \tmsamp{If the conditions 
$(g^{+}_{2})$ or $(g^{-}_{2})$ are true, then (\ref{eq-1})
has a weak solution, i.e., there exists $u\in \H1^1(\Rset^N)$
such that $u$ is a weak solution of (\ref{eq-1}).}
\end{The}

\begin{proof}
As we already know that $\Phi\in C^{1}(X,\Rset)$, We show that
$\Phi$ satisfies the Palais-Smale condition $(PS)$, i.e, give the
sequence $(u_{n})\subset \H1^{1}(\Rset^{N})$ with
$$|\Phi(u_{n})|\leq c ~\mbox{and}~ \Phi^{'}(u_{n})\rightarrow 0$$
we have that $(u_{n})$ has been a subsequence convergent.
Indeed, since
$$\Phi(u_{n})=\frac{1}{2}(Lu_{n},u_{n})_{\H1^{1}(\Rset^{N})}-\int_{\Rset^N}G(x,u_{n})\dd
x,$$
then
$$(\Phi^{'}(u_{n}))v=(Lu_{n},v)_{\H1^{1}(\Rset^{N})}-\int_{\Rset^N}g(x,u_{n})v\dd
x$$
which implies,
\begin{equation}\label{ps-1}
\left|(\Phi^{'}(u_{n}))v\right|=\left|(L{u_{n}},v)-\int_{\Rset^N}g(x,u_{n})v\dd
  x\right|.
\end{equation}
Furthermore,
$$\left|(\Phi^{'}(u_{n}))v\right|\leq\left\|\Phi^{'}(u_{n})\right\|\|v\|$$
where we have
$$\left|(\Phi^{'}(u_{n}))v\right|\leq\left\|v\right\|, ~\forall~ v\in
\H1^{1}({\Rset^N}),$$ for $n$ large enough, since
$$\Phi^{'}(u_{n})\rightarrow 0\Leftrightarrow
\left\|\Phi^{'}(u_{n})\right\|\rightarrow 0$$
thus,
$\left\|\Phi^{'}(u_{n})\right\|\leq 1$. Now, our next step is to show
that 
\begin{equation}\label{limited}
u_{n}=P_{0}u_{n}+P_{-}u_{n}+P_{+}u_{n}
\end{equation}
is bounded. Note that
\begin{enumerate}
\item[$(i)$] If $v=P_{+}u_{n},$ replacing in (\ref{ps-1}) 
follows that
$$
\left|(\Phi^{'}(u_{n}))P_{+}u_{n}\right|=\left|(L{u_{n}},P_{+}u_{n})-
\int_{\Rset^N}g(x,u_{n})P_{+}u_{n}\dd x\right|.$$
Note that using the same reasoning from Lemma \ref{main-6},
\begin{eqnarray*}
(L{u_{n}},P_{+}u_{n})&=&(L(P_{-}u_{n})+L(P_{0}u_{n})+L(P_{+}u_{n}),P_{+}u_{n})\\
&=&(L(P_{+}u_{n}),P_{+}u_{n}).
\end{eqnarray*}
Hence, using the fact that $L$ is defined positive in $X_{+}$, we have
that
\begin{eqnarray*}
(L{u_{n}},P_{+}u_{n})&=&(L(P_{+}u_{n}),P_{+}u_{n})\\
&\geq&\alpha\|P_{+}u_{n}\|^{2}.
\end{eqnarray*}
Furthermore,
\begin{eqnarray*}
\|P_{+}u_{n}\|&\geq&\left|\Phi^{'}(u_{n})P_{+}u_{n}\right|\\
&\geq&|(L{u_{n}},P_{+}u_{n})|-\left|\int_{\Rset^N}g(x,u_{n})P_{+}u_{n}\dd
  x\right|
\end{eqnarray*}
which implies,
$$\|P_{+}u_{n}\|\geq\alpha\|P_{+}u_{n}\|^{2}-\|P_{+}u_{n}\|_{1,Z}.$$
From Theorem \ref{prima-3}, 
we have that
$$\|P_{+}u_{n}\|\geq\alpha\|P_{+}u_{n}\|^{2}-C\|P_{+}u_{n}\|,$$
thus $\|P_{+}u_{n}\|$ is bounded.
\item[$(ii)$] Again, considering $v=P_{-}u_{n}$ in \ref{ps-1}, we obtain
$$(\Phi^{'}(u_{n}))(P_{-}{u_{n}})\leq(Lu_{n},P_{-}u_{n})+\int_{\Rset^N}
|g(x,u_{n})||P_{+}u_{n}|\dd x,$$
from which follow that
$$(\Phi^{'}(u_{n}))(P_{-}u_{n})\leq-\alpha\|P_{-}u_{n}\|^{2}+\|P_{-}u_{n}\|_{1,Z}$$
thus,
\begin{eqnarray*}
-\left\|P_{-}u_{n}\right\|&\leq&\Phi^{'}(u_{n})P_{-}u_{n}\\
&\leq&-\alpha\|P_{-}u_{n}\|^{2}+C_{1}\|P_{-}u_{n}\|
\end{eqnarray*}
consequently,
$$\|P_{-}u_{n}\|\geq\alpha\|P_{-}u_{n}\|^{2}-C_{1}\|P_{-}u_{n}\|.$$
Thus, $\|P_{-}u_{n}\|$ is bounded. 
\end{enumerate}
Now, from $(i)$ and $(ii)$ has been
\begin{equation}\label{ps-2}
\|u_{n}-P_{0}u_{n}\|=\|P_{+}u_{n}+P_{-}u_{n}\|\leq K, ~K\in\Rset_+.
\end{equation}
For another side, we can write
\begin{eqnarray*}
\Phi(u_{n})&=&\frac{1}{2}(L(u_{n}-P_{0}u_{n}),u_{n}-P_{0}u_{n})-\\
&& -\int_{\Rset^{N}}[G(x,u_{n})-G(x,P_{0}u_{n})\dd
x-\int_{\Rset^{N}}G(x,P_{0}u_{n})\dd x
\end{eqnarray*}
we have that,
\begin{eqnarray*}
\int_{\Rset^{N}}G(x,P_{0}u_{n})\dd
x&=&\frac{1}{2}(L(u_{n}-P_{0}u_{n}),u_{n}-P_{0}u_{n})-\Phi(u_{n})-\\
&&-\int_{\Rset^{N}}[G(x,u_{n})-G(x,P_{0}u_{n})]\dd
x.
\end{eqnarray*}
Now, since $L=I-\lambda_{k}S$, we have that $L$ is bounded, thus
\begin{eqnarray*}
|(L(u_{n}-P_{0}u_{n}),u_{n}-P_{0}u_{n})|&\leq&\|L(u_{n}-P_{0}u_{n})\|\|u_{n}-P_{0}u_{n}\|\\
&\leq&\|L\|\|u_{n}-P_{0}u_{n}\|^{2},
\end{eqnarray*}
and from (\ref{ps-2}), we obtain
$$|(L(u_{n}-P_{0}u_{n}),u_{n}-P_{0}u_{n})|\leq C, ~C\in\Rset_+.$$
Note that
\begin{eqnarray*}
\left|\int_{\Rset^{N}}G(x,P_{0}u_{n})\dd
  x\right|&\leq&\frac{1}{2}\left|\left(L(u_{n}-P_{0}u_{n}),u_{n}-P_{0}u_{n}\right)\right|+\left|\Phi(u_{n})\right|+\\
&&+\left|\int_{\Rset^{N}}\left[G(x,u_{n})-G(x,P_{0}u_{n})\right]dx\right|
\end{eqnarray*}
therefore,
$$\left|\int_{\Rset^{N}}G(x,P_{0}u_{n})dx\right|\leq C_{1}, ~C_1\in\Rset_+.$$
Hence, using the condition $(g^{-}_{2})$ we can conclude that 
$\|P_{0}u_{n}\|$ is bounded.
Note that the orthogonality of the projection, we have that
$$\|u_{n}\|^{2}=\|P_{0}u_{n}\|^{2}+\|P_{-}u_{n}\|^{2}+\|P_{+}u_{n}\|^{2},$$
showing that $(u_{n})$ is bounded. We know that
$$
(\Phi^{'}(u))v=\int_{\Rset^{N}}(\nabla u \nabla v +uv)\dd x -
(\psi^{'}(u))v, ~\forall~u,v\in \H1^{1}(\Rset^{N})
$$
where
$$\psi(u)=\int_{\Rset^{N}}\left(\frac{\lambda_{k}}{2}u^{2}h+G(x,u)\right)\dd
x,$$
thus
$$(\nabla(\Phi(u)),v)_{\H1^{1}(\Rset^{N})}=(u,v)_{\H1^{1}(\Rset^{N})}-
(\nabla(\psi(u)),v)_{\H1^{1}(\Rset^{N})}.
$$
Therefore,
$$(\nabla(\Phi(u)),v)_{\H1^{1}(\Rset^{N})}=(u-\nabla(\psi(u)),v)_{\H1^{1}(\Rset^{N})}$$
and consequently,
$\nabla(\Phi(u))=u-\nabla(\psi(u))$. Considering $T(u)=\nabla(\psi(u))$, 
we have that $$\nabla(\Phi(u))=u-T(u).$$ Therefore, 
$$\nabla(\Phi(u_{n}))=u_{n}-T(u_{n})$$ which implies,
$$u_{n}=\nabla(\Phi(u_{n}))+T(u_{n}).$$
Now, since $T:\H1^{1}(\Rset^{N})\to\H1^{1}(\Rset^{N})$
compact (see Appendix C of \cite{MOI}), there exists $(u_{n_{j}})\subset(u_{n})$
such that
$$T(u_{n_{j}})\underset{n_{j}\rightarrow +\infty}{\longrightarrow} u$$
and using the fact that
$$
\Phi^{'}(u_{n})\rightarrow 0\Leftrightarrow \left\|\Phi^{'}(u_{n})\right\|\rightarrow 0 \Leftrightarrow \left\|\nabla\phi(u_{n})\right\|\rightarrow 0\Leftrightarrow \nabla\Phi(u_{n})\rightarrow 0,
$$
passing to the limit in $u_{n_{j}}=\nabla(\Phi(u_{n_{j}}))+T(u_{n_{j}})$
we find
$$u_{n_{j}}\underset{n_{j}\rightarrow\infty}{\longrightarrow} u$$
showing that $\Phi$ satisfies the conditions (PS). 

Finally, we have that $\Phi \in C^{1}(\H1^1(\Rset^{N}),\Rset),$ and
$\Phi$ satisfies the condition of Palais-Smale and using the 
Proposition \ref{main-7} we can apply the Theorem from the
saddle of Rabionowitz with $V=X_{-}$ and $W=X_{0}\oplus X_{+}$
and ensure the existence of a critical point for $\Phi$, therefore, a
weak solution for the problem (\ref{eq-1}).
\end{proof}

\addcontentsline{toc}{chapter}{\bibname}
\bibliographystyle{mybibst}
\bibliography{referencia}

\end{document}